\documentclass[12pt,a4paper]{article}

\usepackage[squaren]{SIunits}
\newcommand{\D}{\mathrm{d}}
\newcommand{\E}{\mathrm{e}}

\newcommand{\nhat}{\hat{\mathbf{n}}}
\newcommand{\Pbar}{\bar{P}}
\newcommand{\rhat}{\hat{\mathbf{r}}}
\newcommand{\bsigma}{\boldsymbol{\sigma}}

\usepackage{algorithm2e}
\usepackage{algorithmicx}
\usepackage{algcompatible}
\usepackage{amsmath}
\usepackage{graphicx}
\usepackage{hyperref}

\begin{document}

\title{Numerical evaluation of the Kirchhoff--Helmholtz integral
  outside a sphere}

\author{M. J. Carley}


\maketitle

\begin{abstract}
  A method is presented for the fast evaluation of the transient
  acoustic field generated outside a spherical surface using surface
  data on the sphere. The method employs Lebedev quadratures, which
  are optimal integration on the sphere, and Lagrange interpolation
  and differentiation in an advanced time algorithm for the evaluation
  of the transient field. Numerical testing demonstrates that the
  approach gives near machine-precision accuracy and a speed-up in
  evaluation time which depends on the order of quadrature rule
  employed but breaks even with direct evaluation at a number of field
  points about~1.15 times the number of surface quadrature nodes,
  making the method an efficient means of evaluating the field
  generated by a large number of sources.
\end{abstract}



\section{Introduction}
\label{sec:intro}

Evaluation of the time-dependent acoustic field outside a source
region is a common task in acoustics. Indeed, the community annoyance
which is often the reason for evaluating the acoustic field can
usually be defined in terms of noise at some distance from an
identified source, such as aircraft operating near a built-up
area. 

In principle, given a time-dependent source distribution, evaluation
of the radiated field is a straightforward summation of the
contribution from the source at each point where it is defined. In
practice, if there are a large number of source points, a situation
which arises when the source is given by a fluid-dynamical
calculation, for example, and the field is required at a large number
of positions, the calculation is extremely demanding of computational
resources. An alternative approach is to evaluate the acoustic
quantities on a surface containing the source and then use these
quantities as a boundary condition for propagation into the region
exterior to the surface. This approach is formally exact and the
question is then how best to implement it numerically for efficient
evaluation of the acoustic field.

For a spherical surface which contains the radiating source, many
approaches are available in the frequency domain, but relatively few
techniques exist for evaluation of a transient signal. One approach is
to compute the evolution of the coefficients of a spherical harmonic
expansion of the field as a function of distance from the sphere
center, using a Laplace transform method. This technique, described as
``teleportation'' in the gravity
literature~\cite{benedict-field-lau13,field-lau15}, has been
considered by a number of
authors~\cite{martin16,martin16a,greengard-hagstrom-jiang14}, with a
particular emphasis on numerically stable formulations for the
evolution of the expansion coefficients. The question of stability and
accuracy are especially important when the field is to be evaluated at
large distance from the source region, where finite difference or
finite element methods are infeasible, and where error can accumulate
in the propagation calculation. 

Another approach which avoids the problems of numerical instability
inherent in the inversion of a Laplace transform is the use of a
surface integral method in the time
domain~\cite{carley-ghorbaniasl16}. This approach has been presented
previously by the author and is developed further here to reduce the
computational requirements in evaluating the field outside a spherical
surface. The method can also be applied to the interior problem, the
evaluation of the acoustic field inside a spherical surface due to
sources outside the surface. This allows the ``transfer'' of the field
radiated from a source region to some other region where the field is
to be evaluated, or the evaluation of a field in a specified domain
due to remote sources. Use of an integral formulation makes the
approach reliable at large distances from the source region, avoiding
the problems of error accumulation when the field must be evolved over
the region between the source and field point.

We offer three potential applications for the method of this
paper. The first is in the calculation of radiation from source
distributions generated by computational methods, such as those which
arise in scattering of transient waves or calculations of turbulent
flow. In the case of scattering calculations, the radiating source
distribution is that on the surface of the scattering body. If the
source is spatially discretized at a resolution proportional to
wavelength, the number of sources required will scale as $f^{2}$ where
$f$ is the maximum frequency to be resolved.  In the case of
three-dimensional turbulent flows, if a sufficient number of points
per wavelength is to be maintained, the number of sources scales as
$f^{3}$ and may be of the order of millions: evaluation and
visualization of the field even over a small region becomes
prohibitive without some acceleration
algorithm~\cite{margnat10,margnat-fortune10,%
  croaker-kessisoglou-kinns-marburg13}. The method of this paper
reduces the computational burden of the calculation to a point where
large scale field evaluation becomes feasible.

The second application is in the use of ffowcs Williams--Hawkings
methods for evaluation of sound using aerodynamic data specified on a
permeable surface. This approach has existed for some
time~\cite{francescantonio97,brentner-farassat98,lyrintzis03} and uses
aerodynamic pressure and velocity on a surface containing the source
where this paper uses acoustic pressure and its normal derivative as
source terms. The choice of surface is arbitrary, which will allow the
method of this paper to be used with a change of source variables.

A final application is the evaluation of transient acoustic fields in
time-domain scattering problems. When such a problem is solved using a
time-domain Boundary Element Method (TDBEM), the scattering surface is
discretized into elements which act as time-varying acoustic
sources. The acoustic pressure generated by these sources, and its
normal derivative, must be calculated in order to generate the
boundary condition for the solution at later time steps. Acceleration
of this calculation is a current research
topic~\cite{takahashi-tanigawa-miyazawa22}, since it is a bottleneck
restricting the speed-up of transient scattering codes. The method of
this paper represents an approach which allows the efficient
evaluation of the acoustic field from groups of sources, similar to
the ``Middleman'' method described by Gumerov and
Duraiswami~\cite[p.172--178]{gumerov-duraiswami04}, and may be useful
to researchers developing improved transient scattering methods.

\section{Radiation from spherical surfaces}
\label{sec:spheres}

The basic method for evaluation of the field outside a surface is the
Kirchhoff--Helmholtz integral~\cite[page~182]{pierce89} which gives
the time-dependent field $p$ at a point $\mathbf{x}$ outside a surface
$S$ in terms of the acoustic pressure $p$ and its normal derivative on
the surface:
\begin{align}
  \label{equ:kirchhoff}
  4\pi p(\mathbf{x}, t)
  &=
  \int_{S} 
  \rhat.\nhat_{1}
  \left(
    \frac{\dot{p}_{1}(\mathbf{x}_{1},\tau)}{Rc}
    +
    \frac{p_{1}(\mathbf{x}_{1},\tau)}{R^{2}}
  \right)
  -
  \frac{1}{R}\frac{\partial p_{1}}{\partial n_{1}}
  \,\D S(\mathbf{x}_{1}),\\
  \mathbf{r} &= \mathbf{x} - \mathbf{x}_{1},
  \,
  R = |\mathbf{r}|,
  \,
  \rhat = \mathbf{r}/R,
  \tau = t - R/c. \nonumber
\end{align}
Here, subscript~1 denotes a variable of integration on $S$, and the
normal derivative of pressure on $S$ is
$\partial p_{1}/\partial n_{1}=\mathbf{n}_{1}.\nabla p_{1}$ with the
normal $\mathbf{n}_{1}$ taken to point out of the surface. Speed of
sound is $c$ and an over-dot denotes differentiation with respect to
time. It will be useful later to use the fact that
$\partial/\partial t\equiv\partial/\partial\tau$. If all sources are
contained inside $S$, the field outside $S$ due to those sources is
given by~\autoref{equ:kirchhoff}. The remainder of this section
describes the main elements used in an efficient method for evaluation
of the integral, which reduces the calculation to a series of matrix
multiplications of the source term on $S$.

We make two further observations. The first is
that~\autoref{equ:kirchhoff} is also valid for the evaluation of the
field inside $S$ generated by sources outside it, with a sign change
on the normal derivative. Secondly, if the method is to be used to
``transfer'' or ``teleport'' the source data from one surface to
another, the normal derivative of the pressure must be found on the
second surface. This requires the integral
\begin{align}
  \label{equ:kirchhoff:normal}
  4\pi \frac{\partial p(\mathbf{x}, t)}{\partial n}
  &=
  \int_{S} 
  \begin{array}[t]{ll}
    \displaystyle
    \left(
      \frac{\dot{p}_{1}}{R^{2}c}
      +
      \frac{p_{1}}{R^{3}}
    \right)
    \nhat.\nhat_{1} 
    - & \\
    \displaystyle    
    \left(
      \frac{\ddot{p}_{1}}{Rc^{2}} +
      3\frac{\dot{p}_{1}}{R^{2}c} +
      3\frac{p_{1}}{R^{3}}
    \right)
    \rhat.\nhat_{1}\rhat.\nhat
    + & \\
    \displaystyle
    \left(
      \frac{1}{R^{2}}\frac{\partial p_{1}}{\partial n_{1}}
      +
      \frac{1}{Rc}\frac{\partial \dot{p}_{1}}{\partial n_{1}}      
    \right)
    \rhat.\nhat
    &
    \,\D S(\mathbf{x}_{1}),
  \end{array}
\end{align}
where $\mathbf{n}$ is the normal at the field point $\mathbf{x}$.

The rest of this section describes the basic elements of the method
and how they are combined to evaluate the acoustic field using surface
data. The appendix contains detailed algorithms which can be used for
implementation of the method. 

\subsection{Advanced time method}
\label{sec:advanced}

The first basic element of the integration algorithm is the advanced
time, or source-time dominant,
method~\cite{casalino03,kessler-wagner04} in which the source retarded
time $\tau$ is specified and the observer reception time $t$ is
calculated. This is particularly useful when the source data are
available at discrete time steps as in the method of this paper. 

\begin{figure}[!htbp]
  \centering
  \includegraphics{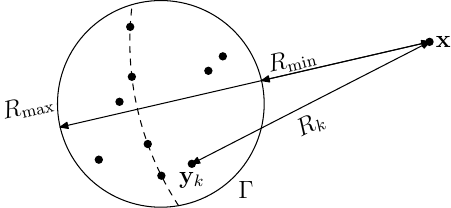}
  \caption{Advanced time evaluation of acoustic field: sources inside
    $S$ have location $\mathbf{y}_{k}$ and sources on the dashed curve
    lie at constant distance from $\mathbf{x}$.}
  \raggedright
\label{fig:advanced}
\end{figure}

\autoref{fig:advanced} shows a set of sources of strength
$q_{k}(\tau)$ at positions $\mathbf{y}_{k}$ inside a boundary
$\Gamma$. The acoustic field at a point $\mathbf{x}$ is
\begin{align}
  \label{equ:advanced}
  p(\mathbf{x},t)
  &=
  \sum_{k}\frac{q_{k}(t-R_{k}/c)}{4\pi R_{k}},
\end{align}
where $R_{k}$ is the distance from $\mathbf{x}$ to the $k$th source
position, see \autoref{fig:advanced}.  The minimum and maximum
distances from $\Gamma$ to $\mathbf{x}$ are $R_{\min}$ and $R_{\max}$
respectively. If $t$ and $\tau$ are discretized with time step
$\Delta t$, the arrival time $t=\tau_{i}+R/c$ of the signal from a
source inside $\Gamma$ can lie anywhere in a range
$\tau_{i}+\Delta R/c\Delta t$, with $\Delta R=R_{\max}-R_{\min}$.

\begin{figure}[!htbp]
  \centering
  \includegraphics{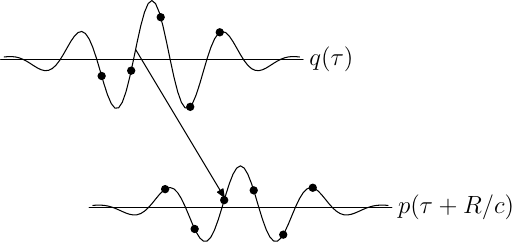}
  \caption{Evaluation of pressure signal by interpolating and scaling
    source term at time $\tau+R/c$}
  \raggedright
\label{fig:advanced:shift}
\end{figure}

\autoref{fig:advanced:shift} shows the principle: a signal generated
at time $\tau_{i}$ contributes to the acoustic pressure at time
$t=\tau_{i}+R/c$ which is not, in general, a time point at which $p$ is
discretized. The contribution to $p$ is thus calculated by
interpolating and scaling $q_{i}=q(\tau_{i})$,
\begin{align}
  p(t) &= \frac{1}{4\pi R}\sum_{k=0}^{K}w_{k}q_{i+k_{0}},
\end{align}
where $w_{k}$ are the weights of a $K+1$ interpolation rule for
evaluation of $q(t-R/c)$. In practice the source data are supplied at
each time step and it is more efficient to accumulate their
contribution to $p$ by incrementing the appropriate elements of
$p$. If $R/c=(n+\delta)\Delta t$, with $0\leq\delta<1$,
\begin{align}
  \label{equ:advanced:increment}
  p_{i+n+k} \gets p_{i+n+k} + w_{k}q(\tau_{i})/4\pi R,\,0\leq k \leq K.
\end{align}
In this paper, we use the Lagrange weights~\cite{berrut-trefethen04}
for interpolation at $x=\delta$ on evenly spaced points with integer
coordinates. When the derivative or second derivative of a source term
contributes to $p$, the corresponding differentiation weights
$\dot{w}_{k}$ and $\ddot{w}_{k}$ are used in place of $w_{k}$,
allowing the evaluation of any terms in the Kirchhoff-Helmholtz
integral, using only the source proper and not its time derivatives,
reducing the memory required for storage of the source terms. The
error in polynomial interpolation of $p$ scales as $(\Delta t)^{K+1}$,
that for $\dot{p}$ as $(\Delta t)^{K}$, and so on. This will have an
effect on the accuracy of evaluation of the signal, especially when
the normal derivative $\partial p/\partial n$ is to be computed, as it
requires evaluation of a second derivative with error proportional to
$(\Delta t)^{K-1}$, requiring the use of a high order interpolation
scheme. 

When multiple sources contribute to the radiated field, contributions
from the same value of $\tau$ will not necessarily contribute to the
same values of $p_{i}$ because of variations in $R/c$. Given a vector
$\bsigma_{i}$ of $n_{s}$ source strengths at retarded time $\tau_{i}$,
the incrementing of $p$ is implemented as a matrix multiplication
\begin{align}
  \label{equ:advanced:matrix}
  p_{i+n+k} \gets p_{i+n+k} + \mathbf{W}\bsigma_{i}, \,0\leq k \leq K,
\end{align}
where $\mathbf{W}$ is an $n_{s}\times(\Delta R/c\Delta t)$ matrix with
each row given by \autoref{equ:advanced:increment}, with zero padding
to ensure the correct alignment of signals from sources at different
distances $R$. Corresponding matrices $\dot{\mathbf{W}}$ and
$\ddot{\mathbf{W}}$ are used to evaluate terms involving time
derivatives.

Finally, we note that sources which lie at the same distance from
$\mathbf{x}$ can have their contributions summed and be treated as a
single source when incrementing $p$. Integration over a spherical
surface is implemented using such a summation.

\subsection{Interpolation on spherical surfaces}
\label{sec:interpolation}

Efficient evaluation of the Kirchhoff--Helmholtz integral depends on a
suitable choice of an interpolation scheme on the sphere. In previous
work~\cite{carley-ghorbaniasl16}, we used spherical harmonic
interpolation based on trapezoidal rule quadrature in azimuth and a
Gaussian quadrature in elevation. In this paper, we again adopt
spherical harmonics as our interpolation functions, but employ Lebedev
quadrature rules~\cite{lebedev77} for the interpolation nodes. These
rules appear to be optimal for integration on the
sphere~\cite{beentjes15} and give a considerable reduction in the
number of surface points where source quantities must be evaluated,
with a corresponding improvement in efficiency and memory use.

Taking the origin of coordinates at the center of the sphere, a point
$\mathbf{x}$ is given by
\begin{align*}
  \mathbf{x} &=
  \rho(\sin\theta\cos\phi,\,\sin\theta\sin\phi,\,\cos\theta), 
\end{align*}
with $\rho=a$ the radius of the sphere. A function on the sphere can
be expressed
\begin{align}
  \label{equ:analysis:sph}
  f(\theta,\phi) 
  &=
  \sum_{n=0}^{\infty}
  \sum_{m=0}^{n}
  \Pbar_{n}^{m}(\cos\theta)
  [a_{n,m}\cos (m\phi) + b_{n,m}\sin (m\phi)],
\end{align}
where $\Pbar_{n}^{m}(\theta)$ is the normalized associated Legendre
function,
\begin{align}
  \label{equ:analysis:pbar}
  \Pbar_{n}^{m}(\theta)
  &=
  \left[
    \frac{2n+1}{2}
    \frac{(n-m)!}{(n+m)!}
  \right]^{1/2}
  P_{n}^{m}(\cos\theta).
\end{align}
The coefficients $a_{n,m}$ and $b_{n,m}$ are given by integration over
the spherical surface, exploiting the orthogonality of the spherical
harmonics: 
\begin{subequations}
  \label{equ:analysis:cfft}
  \begin{align}
    a_{n,m} 
    &=
    \frac{1}{1+\delta_{m0}}
    \int_{0}^{\pi}\int_{0}^{2\pi}
    \Pbar_{n}^{m}(\cos\theta)
    \cos (m\phi)
    f(\theta,\phi)
    \,\D\phi\,\sin\theta\,\D\theta,\\
    b_{n,m}
    &=
    \int_{0}^{\pi}\int_{0}^{2\pi}
    \Pbar_{n}^{m}(\cos\theta)
    \sin (m\phi)
    f(\theta,\phi)
    \,\D\phi\,\sin\theta\,\D\theta,
  \end{align}
\end{subequations}
where $\delta_{ij}$ is the Kronecker delta. 

Adopting the Lebedev rules~\cite{lebedev77}, which are symmetric and
integrate spherical polynomials exactly up to some specified order
$N$, the expansion coefficients are given by
\begin{subequations}
  \label{equ:interpolation:expansion}
  \begin{align}
    a_{n,m} 
    &=
    \frac{1}{1+\delta_{m0}}
    \sum_{i=1}^{N_{Q}}
    w_{i}\Pbar_{n}^{m}(\cos\theta_{i})\cos (m\phi_{i}) f(\theta_{i},\phi_{i}),\\
    b_{n,m} 
    &=
    \sum_{i=1}^{N_{Q}}
    w_{i}\Pbar_{n}^{m}(\cos\theta_{i})\sin (m\phi_{i}) f(\theta_{i},\phi_{i}),
\end{align}
\end{subequations}
where $N_{Q}$ is the number of nodes in the Lebedev rule.  The
evaluation is implemented as a matrix multiplication,
\begin{align}
  \label{equ:interpolation:matrix}
  \mathbf{a} &= \mathbf{A}\mathbf{f},
\end{align}
where the elements of matrix $\mathbf{A}$ are given by
\autoref{equ:interpolation:expansion}, and the vector $\mathbf{f}$
holds the values of the function to be interpolated at the quadrature
nodes. Algorithm~\ref{alg:interpolation:matrix} gives details of the
evaluation of $\mathbf{A}$.

To evaluate the interpolant at some point $(\theta,\phi)$
\begin{align}
  \label{equ:interpolation:eval}
  f(\theta,\phi) &\approx \mathbf{b}(\theta,\phi).\mathbf{f},    
\end{align}
where the weight vector $\mathbf{b}$ is found using the spherical
harmonics evaluated at $(\theta,\phi)$:
\begin{align}
  \label{equ:interpolation:b}
  \mathbf{b}(\theta,\phi)
  &=
  [\ldots P_{n}^{m}(\cos\theta)\cos m\phi\quad
  P_{n}^{m}(\cos\theta)\sin m\phi\, \ldots]\mathbf{A},
\end{align}
described in Algorithm~\ref{alg:interpolation}. 

\subsection{Kirchhoff--Helmholtz integral on a sphere}
\label{sec:kirchhoff}

\begin{figure}[!htbp]
  \centering
  \includegraphics{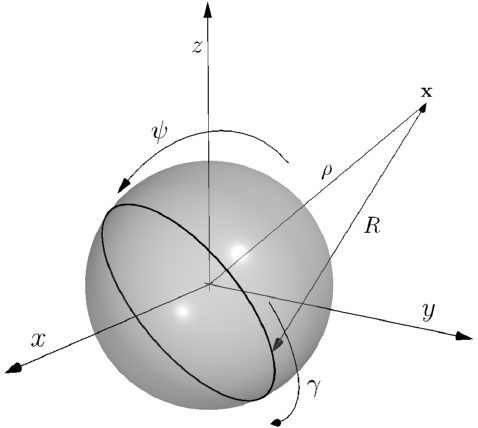}
  \caption{Coordinate system for radiating surface and rotated system
    for integral evaluation}
  \raggedright
\label{fig:rotation}
\end{figure}

Integration over a spherical surface is performed in a spherical polar
coordinate system aligned with the vector to the field point
$\mathbf{x}$, \autoref{fig:rotation}, such that
$\mathbf{x}=(0,0,\rho)$. Points with elevation $\psi$ in the new
coordinate system lie at constant distance $R$ from $\mathbf{x}$ and
an azimuthal angle $\gamma$ completes the specification of
position. In this coordinate system, the Kirchhoff--Helmholtz integral
becomes
\begin{align}
  \label{equ:kirchhoff:spherical}
  p(\mathbf{x}, t)
  &=
  \frac{a^{2}}{4\pi}\frac{\partial}{\partial t}
  \int_{0}^{\pi}
  \frac{f_{0}}{Rc}P_{0}
  \sin\psi
  \,\D\psi + \frac{a^{2}}{4\pi}
  \int_{0}^{\pi}
  \frac{f_{0}}{R^{2}}P_{0}
  \sin\psi
  \,\D\psi\nonumber\\
  &-
  \frac{a^{2}}{4\pi}
  \int_{0}^{\pi}
    \frac{1}{R} N_{0}
  \sin\psi
  \,\D\psi,\\
  R^{2} &= \rho^{2} + a^{2} - 2\rho a\cos\psi,\nonumber
\end{align}
and if a normal derivative is required at some point,
\begin{align}
  \label{equ:kirchhoff:normal:1}
  \frac{\partial p(\mathbf{x}, t)}{\partial n}
  &=
  \frac{a^{2}}{4\pi}\frac{\partial^{2}}{\partial t^{2}}
  \int_{0}^{\pi}
  \frac{f_{0}}{Rc^{2}}
  \left(
    f_{1}P_{C} + 
    f_{2}P_{S} +
    f_{3}P_{0}      
  \right)
  \sin\psi \,\D\psi\nonumber\\
  &+ 
  \frac{a^{2}}{4\pi}\frac{\partial}{\partial t}
  \int_{0}^{\pi}
  \frac{1}{R^{2}c}
  \left(
    n_{x}\sin\psi P_{C} + n_{y}\sin\psi P_{S} + n_{z}\cos\psi P_{0}
  \right)
  \sin\psi\,\D\psi\nonumber\\
  &+
  \frac{a^{2}}{4\pi}\frac{\partial}{\partial t}
  \int_{0}^{\pi}
  \frac{3f_{0}}{R^{2}c}
  \left(
    f_{1}P_{C} +
    f_{2}P_{S} +
    f_{3}P_{0}
  \right)
  \sin\psi\,\D\psi\nonumber\\
  &+
  \frac{a^{2}}{4\pi}
  \int_{0}^{\pi}
  \frac{1}{R^{3}}
  \left(
    n_{x}\sin\psi P_{C} + n_{y}\sin\psi P_{S} + n_{z}\cos\psi P_{0}
  \right)
  \sin\psi\,\D\psi\nonumber\\
  &+
  \frac{a^{2}}{4\pi}
  \int_{0}^{\pi}
  \frac{3f_{0}}{R^{3}}
  \left(
    f_{1}P_{C} +
    f_{2}P_{S} +
    f_{3}P_{0}
  \right)
  \sin\psi\,\D\psi\nonumber\\
  &-
  \frac{a^{2}}{4\pi}
  \frac{\partial}{\partial t}
  \int_{0}^{\pi}
  \frac{1}{Rc}
  \left(f_{1}N_{C} + f_{2}N_{S} + f_{3}N_{0}\right)
  \sin\psi \,\D\psi\nonumber\\
  &-
  \frac{a^{2}}{4\pi}
  \int_{0}^{\pi}
  \frac{1}{R^{2}}
  \left(f_{1}N_{C} + f_{2}N_{S} + f_{3}N_{0}\right)
  \sin\psi \,\D\psi,
\end{align}
with the normal $\mathbf{n}=(n_{x},n_{y},n_{z})$ given in the rotated
coordinate system. The intermediate quantities are integrals of source
terms over $\gamma$
\begin{subequations}
  \label{equ:azimuthal}
\begin{align}
  \left[
    \begin{array}{l}
      P_{0}(\psi, \tau)\\
      P_{C}(\psi, \tau)\\
      P_{S}(\psi, \tau)      
    \end{array}
  \right]
  &=
  \int_{0}^{2\pi}
  p_{1}(\psi, \gamma, \tau)
  \left[
    \begin{array}{l}
      1\\
      \cos\gamma\\
      \sin\gamma
    \end{array}
  \right]
  \,\D\gamma,\\
  \left[
    \begin{array}{l}
      N_{0}(\psi, \tau)\\
      N_{C}(\psi, \tau)\\
      N_{S}(\psi, \tau)      
    \end{array}
  \right]
  &=
  \int_{0}^{2\pi}
  \frac{\partial}{\partial n_{1}}p_{1}(\psi, \gamma, \tau)
  \left[
    \begin{array}{l}
      1\\
      \cos\gamma\\
      \sin\gamma
    \end{array}
  \right]
  \,\D\gamma,\\  
  f_{0} &= \frac{\rho\cos\psi - a}{R},\,
  f_{1} = n_{x}\frac{a\sin\psi}{R},\,
  f_{2} = n_{y}\frac{a\sin\psi}{R},\,
  f_{3} = n_{z}\frac{a\cos\psi - \rho}{R}.
\end{align}  
\end{subequations}

Integration is performed using a Gauss--Legendre quadrature in $\psi$
and a trapezoidal rule in $\gamma$. The trapezoidal rule is
implemented at a given $\psi_{i}$ as a scalar product with appropriate
weight vectors:
\begin{align}
  \int_{0}^{2\pi}f(\theta,\phi)
  \left[
    \begin{array}{l}
      1\\
      \cos\gamma\\
      \sin\gamma
    \end{array}
  \right]
  \,\D\gamma
  &\approx
  \left[
    \begin{array}{l}
      \mathbf{q}.\mathbf{f}\\
      \mathbf{q}_{C}.\mathbf{f}\\
      \mathbf{q}_{S}.\mathbf{f}
    \end{array}
  \right],\\
  \label{equ:integral:B}
  \left[
    \begin{array}{l}
      \mathbf{q}\\
      \mathbf{q}_{C}\\
      \mathbf{q}_{S}
    \end{array}
  \right]
  &=
  \frac{2\pi}{N_{\gamma}}\sum_{j=0}^{N_{\gamma}-1}
  \mathbf{b}(\theta(\psi_{i},\gamma_{j}),\phi(\psi_{i},\gamma_{j})
  \left[
    \begin{array}{l}
      1\\
      \cos\gamma\\
      \sin\gamma
    \end{array}
  \right],\quad
  \gamma_{j} = 2\pi (j-1)/N_{\gamma}.
\end{align}
Algorithm~\ref{alg:integration} details the generation of the vectors
$\mathbf{q}$, $\mathbf{q}_{C}$, and $\mathbf{q}_{S}$. 

We note that for the case of a field point on the surface of the
sphere, $\rho=a$, a hypersingular quadrature rule can be
used~\cite{paget81} to deal with the singular integrand in
\autoref{equ:kirchhoff:normal:1}.

\subsection{Implementation}
\label{sec:implementation}

The basic elements of the previous sections can now be combined into a
method for the evaluation of the field, and its normal derivative if
required, radiated to a point $\mathbf{x}$ from a sphere of radius $a$
centered at the origin. At time step $i$, vectors $\bsigma_{i}$ and
$\bsigma^{(n)}_{i}$ contain the surface acoustic pressure and its
gradient respectively. Each vector is of length $N_{Q}$ where $N_{Q}$
is the number of surface interpolation nodes from the Lebedev
quadrature.

Given a quadrature rule of length $N_{\psi}$ for integration over
$\psi$, the contribution of $\bsigma$ at any time step to the radiated
field can be evaluated by a matrix multiplication which encodes the
integrations of \autoref{equ:kirchhoff:spherical}. The matrices are
products of the matrices which encode the advanced time calculations
of \autoref{sec:advanced} and the azimuthal integrations of
\autoref{sec:kirchhoff}. Details of the intermediate stages and of the
resulting matrices are given in Algorithms~\ref{alg:matrices}
and~\ref{alg:shifts} in the appendix. The method for incrementing the
radiated field at each time step is then
\begin{align}
  \label{equ:implementation:matrix}
  p_{i+n+k} \gets p_{i+n+k} + \mathbf{S}\bsigma_{i} +
  \mathbf{S}_{N}\bsigma_{i}^{(n)},\\
  \left(\frac{\partial p}{\partial n}\right)_{i+n+k} \gets
  \left(\frac{\partial p}{\partial n}\right)_{i+n+k} +
  \overline{\mathbf{S}}\bsigma_{i} + 
  \overline{\mathbf{S}}_{N}\bsigma_{i}^{(n)}.
\end{align}
The matrices $\mathbf{S}$, $\mathbf{S}_{N}$, etc. are approximately of
size $(K+1)\times N_{Q}$, where the exact number of rows depends on
the difference in transit time between the field point and the
furthest and closest points on the spherical surface. In any case, the
size of the matrices is independent of the number of quadrature nodes
in $\psi$ and $\gamma$ and the precision of the interpolation has an
effect only insofar as it sets the length of Lebedev quadrature rule
required to accurately integrate spherical harmonics up to order
$N_{S}$. 



\subsection{Computational demands}
\label{sec:demands}

The method presented in this paper is a development of previous
work~\cite{carley-ghorbaniasl16}. Improvements are the use of Lebedev
quadratures in place of tensor-product integration on the sphere and
the use of Lagrange interpolation to evaluate time derivatives rather
than storing the value of $\dot{p}$. Beentjes~\cite{beentjes15} gives
a useful comparison of the computational efficiency of different
quadrature rules on the sphere, which allows an estimate of the memory
and computation demands of the two approaches. For the previous
method, using Gauss-Legendre quadrature in $\theta$ and a trapezoidal
rule in $\phi$, the quadrature rule requires $(N_{s}+1)^{2}/2$ points
to resolve spherical harmonics up to order $N_{S}$. With three source
terms, $p$, $\partial p/\partial n$, and $\dot{p}$, the memory
requirement per time step is $3(N_{S}+1)^{2}/2$. For the higher order
Lebedev quadrature rules, which are near optimal, the number of
quadrature nodes is approximately $(N_{S}+1)^{2}/3$. With only two
source terms required, $p$ and $\partial p/\partial n$, the memory
requirement is $2(N_{S}+1)^{2}/3$, less than half that required for
the previous method. While the storage required per time step is
relatively modest, halving the requirement doubles the number of time
steps which can be held in memory on a given machine. 

The computational demand is dominated by the cost of the matrix
multiplications of \autoref{equ:implementation:matrix}. At fixed time
interpolation order $K$, the cost of matrix multiplication scales as
the number of data nodes squared. The cost per time step of the matrix
multiplications for three source terms in the original method is then
proportional to $3(N_{S}+1)^{4}/4$; the corresponding cost for the two
matrix multiplications in the method of this paper is
$2(N_{S}+1)^{4}/9$, less than one third of the cost in the original
algorithm. We note that this reduction in memory and computation time
is achieved with no loss of accuracy.

\section{Results}
\label{sec:results}

\begin{figure}[!htbp]
  \centering
  \includegraphics{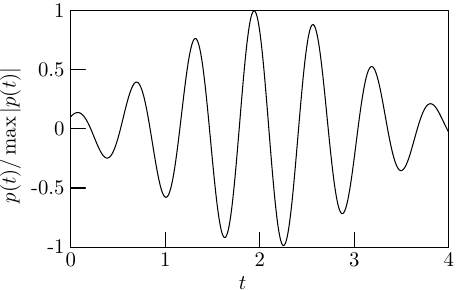}
  \caption{Normalized signal from single source with strength given
    by~\autoref{equ:results:source}.}
  \raggedright
\label{fig:results:signal}
\end{figure}

The accuracy and the speed of the evaluation method depend on a number
of parameters. For concision, we present results for two test
cases. In the first case, where we examine the accuracy and
convergence of the technique, the surface pressure data are generated
using a single point source placed inside a spherical surface of
radius $a=1$. The source position is $0.7\times(1,-1,1)/\sqrt{3}$ and
its strength is
\begin{align}
  \label{equ:results:source}
  q(\tau) &= \E^{-\alpha(\tau-t_{0})^{2}}\cos\Omega\tau,
\end{align}
with $\alpha=1/2$, $t_{0}=2$, $\Omega=10$, see
\autoref{fig:results:signal}. The radiated field $p$ is evaluated at a
radius $\rho=2$ with varying time interpolation order, and varying
number of spherical harmonics in the interpolation scheme, for
$0\leq\tau\leq4$, with varying number of time steps $n_{t}$. The
normal derivative $\partial p/\partial n$ is evaluated at the same
point, with an arbitrarily chosen normal. Error in the computed signal
is
\begin{align}
  \label{equ:results:error}
  \epsilon
  &=
  \frac{\max|p_{d}(t)-p_{c}(t)|}{\max|p_{d}(t)|},
\end{align}
where subscript $d$ denotes the exact pressure evaluated directly from
the source data and $c$ that computed using the surface integral
method.

\begin{figure}[!htbp]
  \centering
  \includegraphics{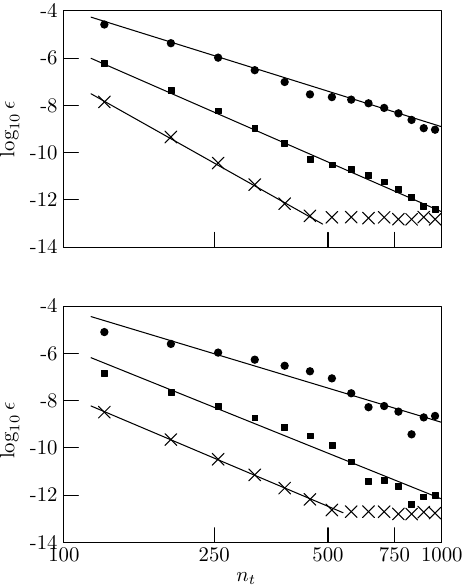}
  \caption{Error $\epsilon$ in $p$ (upper plot) and
    $\partial p/\partial n$ (lower plot) for single source input
    against number of time steps $n_{t}$ with time interpolation
    order~4 (bullets),~6 (squares), and~8 (crosses); straight lines
    have slope~$-5$,~$-7$,~$-9$ on log axes in upper plot
    and~$-4.8$,~$-6.4$,~$-6.8$ on log axes in lower plot,
    respectively.}
  \raggedright
\label{fig:results:time}
\end{figure}

\autoref{fig:results:time} shows error in the computed time record for
$p$ and $\partial p/\partial n$ as a function of number of time steps
and temporal interpolation order. The upper plot, of the error in $p$
shows steady convergence as the time step is reduced, with the higher
order interpolation schemes reaching a relative error of
about~$10^{-13}$ and, not shown here, absolute error of about machine
precision. The error in $p$ reduces at a rate approximately equal to
$(\Delta t)^{K+1}$, where $K$ is the order of the Lagrange
interpolant, in line with the error behavior for polynomial
interpolation. The error in $\partial p/\partial n$ reduces more
slowly probably because of the second derivative which is evaluated to
a lower order of accuracy than the lower order terms. The error in
$\partial p/\partial n$ for $K=8$ reduces more slowly than might be
expected, though it still converges quite rapidly. This may be caused
by some instability in the high order polynomial differentiation on
evenly spaced points, the so-called Runge
phenomenon~\cite{berrut-trefethen04}.

\begin{figure}[!htbp]
  \centering
  \includegraphics{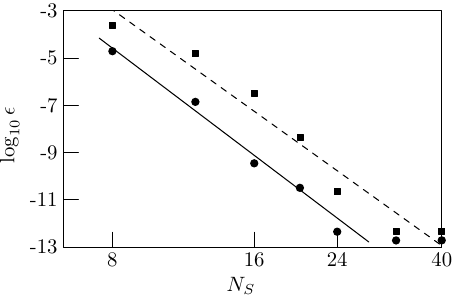}
  \caption{Error $\epsilon$ for single source input. Upper plot: error
    against maximum order of spherical harmonics for evaluation of $p$
    (circles) and $\partial p/\partial n$ (squares); straight lines
    have slope~-15 and~-14.3 respectively on log axes. Lower plot:
    error against field point radius; circles $N_{S}=8$, squares
    $N_{S}=12$, crosses $N_{S}=16$, diamonds $N_{S}=20$, upward
    triangles $N_{S}=24$, downward triangles $N_{S}=32$.}
  \raggedright
\label{fig:results:nsph}
\end{figure}

\autoref{fig:results:nsph} shows the error in $p$ and
$\partial p/\partial n$ as a function of the maximum order of
spherical harmonics in the interpolation at a fixed radius $\rho=2a$,
and as expected there is rapid convergence with machine precision
being achieved when the order of the spherical harmonic expansion is
about~32.

The lower plot in \autoref{fig:results:nsph} shows the error in $p$ as
a function of distance from the source sphere with varying order of
spherical harmonic expansion and it can be seen to decay very rapidly,
even for a relatively low order interpolation scheme, $N_{S}=8$.

To assess the computational effort required in applying the field
evaluation scheme, we report the number of field points at which the
scheme breaks even with direct evaluation as a function of the number
of sources $n_{s}$. If direct evaluation of the field at one point
takes time $t_{d}$ and surface integration takes $t_{s}$ with a
pre-processing time $t_{p}$, the break-even number of field points is
\begin{align}
  n_{f} &= \frac{t_{p}}{t_{d}-t_{s}}.
\end{align}

\begin{figure}[!htbp]
  \centering
  \includegraphics{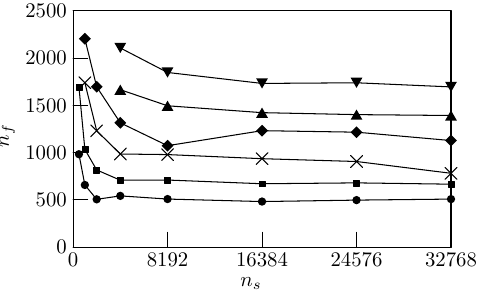}
  \caption{Performance of the method as a function of number of
    sources, $n_{s}$. Upper plot: break-even number of field points
    against number of sources for varying number of sphere quadrature
    points:~434 (bullets),~590 (squares),~770 (crosses),~974
    (diamonds),~1202 (upward triangles),~1454 (downward
    triangles). Lower plot: error in $p$ evaluated at field point
    radius~2, symbols as in upper plot.}
  \raggedright
\label{fig:results:multi}
\end{figure}

The second test case, to assess computation time and accuracy as a
function of interior source number, uses randomized sources placed at
random positions inside a ball of radius~0.7 and a spherical surface
of radius $a=1$,
\begin{align*}
  p(\mathbf{x}) &=
  \sum_{k=1}^{n_{s}}
  \frac{q_{k}(t-R_{k}/c)}{4\pi R_{k}},\\
  q_{k} &= \cos \Omega_{k}\tau,\quad 9\leq\Omega_{k}\leq 11,\\
  R_{k} &= |\mathbf{x}-\mathbf{y}_{k}|,
\end{align*}
with $\Omega_{k}$ and $\mathbf{y}_{k}$ randomly assigned.

Calculations are run using a fixed order of spherical harmonic
interpolation, $N_{S}=32$, and varying orders of Lebedev quadrature,
with the number of sphere quadrature nodes varying from~434 to~1454,
measuring $t_{d}$, $t_{s}$ and $t_{p}$. The upper plot of
\autoref{fig:results:multi} shows the break-even number of field
points $n_{f}$ as a function of $n_{s}$. At large $n_{s}$, the
break-even value of $n_{f}$ becomes roughly constant at a value of
about~1.15 times the number of quadrature points on the surface,
significantly less than the number of sources, despite the high
accuracy of the method, shown in the lower plot of
\autoref{fig:results:multi} which gives the error in $p$ at a field
radius of~2 as a function of number of quadrature points. As in the
single-source case, see \autoref{fig:results:nsph}, the error reduces
rapidly with distance from the spherical surface. 

\section{Conclusions}
\label{sec:conclusions}

A method has been presented for the evaluation of the acoustic field
outside a source region, based on an efficient technique for the
evaluation of the Kirchhoff--Helmholtz integral on a spherical
surface. Compared to an earlier version of the method, the new
approach significantly reduces the computational time and memory
required, while maintaining high accuracy. This is achieved through
the use of high order schemes for temporal interpolation and
differentiation and the adoption of efficient quadrature rules for
interpolation on the sphere. Calculation of the number of field
evaluations at which the method breaks even with direct evaluation
shows that even for relatively modest source numbers, in the low
thousands, the radiated field can be accurately evaluated very much
faster using the method of this paper.


\appendix*

\section{Algorithms}
\label{sec:algorithms}

For convenience, we present a summary of the elements of the algorithm
in a form suitable for implementation. Vectors are represented
throughout as lower-case bold letters, and matrices by upper-case bold
letters. Details are given for those elements of the method which are
required for the evaluation of both $p$ and $\partial p/\partial
n$. If only the acoustic pressure is required, the quantities
appropriate to the normal derivative of pressure may be neglected.


\begin{algorithm}
  \caption{Matrix for evaluation of spherical harmonic interpolation
    coefficients }
  \label{alg:interpolation:matrix}
  \begin{algorithmic}
    \STATE input Lebedev quadrature rule of $N_{Q}$ points, maximum
    order of spherical harmonics $N_{S}$
    \FOR{$i\gets1,\ldots,N_{Q}$}
    \FOR{$n\gets 0, \ldots, N_{S}$}
    \FOR{$m\gets 0,\ldots,n$}
    \STATE $j\gets n(n+1)/2+m$
    \STATE $A_{2j+1,i}\gets w_{i}\Pbar_{n}^{m}(\theta_{i})\cos m\phi_{i}$
    \STATE $A_{2j+2,i}\gets w_{i}\Pbar_{n}^{m}(\theta_{i})\sin m\phi_{i}$
    \ENDFOR
    \ENDFOR    
    \ENDFOR
    \STATE output matrix $\mathbf{A}$ of size
    $(N_{S}+1)(N_{S}+2)\times N_{Q}$ which evaluates coefficients of
    spherical harmonic expansion
  \end{algorithmic}
\end{algorithm}

\begin{algorithm}
  \caption{Vector for interpolation of function on sphere at
    $(\theta,\phi)$}
  \label{alg:interpolation}
  \begin{algorithmic}
    \STATE set $\theta$, $\phi$
    \STATE initialize $\mathbf{b}$
    \FOR{$n\gets 0, \ldots, N_{S}$}
    \FOR{$m\gets 0,\ldots,n$}
    \STATE $i\gets n(n+1)/2+m$
    \STATE $b_{2i+1}\gets \Pbar_{n}^{m}(\theta)\cos m\phi$
    \STATE $b_{2i+2}\gets \Pbar_{n}^{m}(\theta)\sin m\phi$
    \ENDFOR
    \ENDFOR
    \STATE output vector $\mathbf{b}$ of length $(N_{S}+1)(N_{S}+2)$
    containing spherical harmonics evaluated at $(\theta,\phi)$ for
    use with interpolation coefficients
  \end{algorithmic}
\end{algorithm}

\begin{algorithm}
  \caption{Weights for integration in $\gamma$ at given $\psi$}
  \label{alg:integration}
  \begin{algorithmic}
    \STATE input field point coordinates $(\theta, \phi)$, $\psi$, number
    of quadrature nodes in $\gamma$, $N_{\gamma}$, matrix $\mathbf{A}$
    from Algorithm~\ref{alg:interpolation:matrix}
    \STATE set temporary vectors $\mathbf{q}$, $\mathbf{q}_{C}$,
    $\mathbf{q}_{S}$ to zero
    \FOR{$i\gets0,N_{\gamma}-1$}
    \STATE $\gamma_{i}\gets 2\pi i/N_{\gamma}$
    \STATE find transformed coordinates $\theta_{i}(\psi,\gamma_{i})$,
    $\phi_{i}(\psi,\gamma_{i})$,
    \STATE generate interpolation vector $\mathbf{b}$ using
      Algorithm~\autoref{alg:interpolation}
    \STATE $\mathbf{q}\gets \mathbf{q}+\mathbf{b}$
    \STATE $\mathbf{q}_{C}\gets \mathbf{q}_{C}+\mathbf{b}\cos\gamma_{i}$
    \STATE $\mathbf{q}_{S}\gets \mathbf{q}_{S}+\mathbf{b}\sin\gamma_{i}$
    \ENDFOR
    \STATE \COMMENT{vectors $\mathbf{q}$, $\mathbf{q}_{C}$,
      $\mathbf{q}_{S}$ perform integration in $\gamma$ based on
      interpolation coefficients: now convert to integration based on
      nodal values}
    \STATE $\mathbf{q}\gets \mathbf{A}^{T}\mathbf{q}2\pi/N_{\gamma}$
    \STATE $\mathbf{q}_{C}\gets \mathbf{A}^{T}\mathbf{q}_{C}2\pi/N_{\gamma}$
    \STATE $\mathbf{q}_{S}\gets \mathbf{A}^{T}\mathbf{q}_{S}2\pi/N_{\gamma}$
    \STATE output vectors $\mathbf{q}$, $\mathbf{q}_{C}$,
    $\mathbf{q}_{S}$ of length $N_{Q}$ which perform the integrations
    of \autoref{equ:azimuthal}
  \end{algorithmic}
\end{algorithm}

\begin{algorithm}[ht]
  \caption{Matrices for advanced time field evaluation}
  \label{alg:matrices}
  \begin{algorithmic}
    \STATE input sphere radius $a$, field point distance $\rho$, maximum
    order of spherical harmonics $N_{S}$, order of time interpolation
    $K$, length of quadrature rule $N_{\psi}$, speed of sound $c$,
    normal $\mathbf{n}=(n_{x},n_{y},n_{z})$
    \FOR{$i=1,\ldots,N_{\psi}$}
    \STATE set quadrature node $\psi_{i}$
    \STATE $R_{i}=(\rho^{2}+a^{2} - 2a\rho\cos\psi_{i})^{1/2}$
    \STATE $N_{i}=\lfloor R_{i}/c/\Delta t\rfloor$,
    $\delta_{i}=(R_{i}/c-N_{i}\Delta t)/\Delta t$
    \STATE evaluate $w_{k}$, $\dot{w}_{k}$, $\ddot{w}_{k}$, $k=0,\ldots,K$
    \COMMENT{Lagrange interpolation and differentiation weights using
      method of Berrut and Trefethen~\cite{berrut-trefethen04} for 
      $x=1-\delta$ and reversing the order of the weights}
    \STATE $\dot{w}_{k}\gets \dot{w_{k}}/\Delta t$
    \STATE $\ddot{w}_{k}\gets \ddot{w_{k}}/(\Delta t)^{2}$
    \STATE $c_{0}\gets (\rho\cos\psi_{i}-a)/R_{i}$,
    $c_{1}\gets an_{x}\sin\psi_{i}/R_{i}$,
    $c_{2}\gets an_{y}\sin\psi_{i}/R_{i}$,
    $c_{3}\gets n_{z}(a\cos\psi_{i}-\rho)/R_{i}$
    \FOR{$k=0,\ldots,K$}
    \STATE $j \gets N_{i}-N_{1}+k+1$
    \STATE $W_{ji} \gets c_{0}\dot{w}_{k}/Rc +
    c_{0}w_{k}/R_{i}^{2}$
    \STATE $V_{ji} \gets -w_{k}/R$
    \STATE $\overline{W}^{(0)}_{ji}\gets
    c_{0}c_{3}\ddot{w}_{k}/R_{i}c^{2} +
    (n_{z}\cos\psi_{i} + 3c_{0}c_{3})\dot{w}_{k}/R_{i}^{2}c +
    (n_{z}\cos\psi_{i} + 3c_{0}c_{3})\dot{w}_{k}/R_{i}^{3}
    $
    \STATE $\overline{W}^{(C)}_{ji}\gets
    c_{0}c_{1}\ddot{w}_{k}/R_{i}c^{2} +
    (n_{x}\sin\psi_{i} + 3c_{0}c_{1})\dot{w}_{k}/R_{i}^{2}c +
    (n_{x}\sin\psi_{i} + 3c_{0}c_{1})\dot{w}_{k}/R_{i}^{3}
    $
    \STATE $\overline{W}^{(S)}_{ji}\gets
    c_{0}c_{2}\ddot{w}_{k}/R_{i}c^{2} +
    (n_{y}\sin\psi_{i} + 3c_{0}c_{2})\dot{w}_{k}/R_{i}^{2}c +
    (n_{y}\sin\psi_{i} + 3c_{0}c_{2})\dot{w}_{k}/R_{i}^{3}
    $
    \STATE $\overline{V}^{(0)}_{ji}=-c_{3}\dot{w}_{k}/Rc - c_{3}w_{k}/R_{i}^{3}$
    \STATE $\overline{V}^{(C)}_{ji}=-c_{1}\dot{w}_{k}/Rc - c_{1}w_{k}/R_{i}^{3}$
    \STATE $\overline{V}^{(S)}_{ji}=-c_{2}\dot{w}_{k}/Rc - c_{2}w_{k}/R_{i}^{3}$
    \ENDFOR
    \ENDFOR
    \STATE output $N_{W}\times N_{\psi}$
    matrices $\mathbf{W}$, $\mathbf{V}$, $\overline{\mathbf{W}}^{(0,C,S)}$,
    $\overline{\mathbf{V}}^{(0,C,S)}$, with $N_{W}=(N_{N_{\psi}}-N_{1}+K+1)$
  \end{algorithmic}
\end{algorithm}

\begin{algorithm}
  \caption{Matrices to increment radiated field from source data}
  \label{alg:shifts}
  \begin{algorithmic}
    \STATE input $a$, $c$, $(\rho,\theta,\phi)$, matrices $\mathbf{W}$,
    $\mathbf{N}$, $\overline{\mathbf{W}}^{(0,C,S)}$,
    $\overline{\mathbf{V}}^{(0,C,S)}$, from
    Algorithm~\ref{alg:matrices}, Gauss-Legendre quadrature rule
    $(\psi_{i},w_{i}^{(\psi)})$ of length $N_{\psi}$, length of
    trapezoidal rule $N_{\gamma}$ 
    \STATE initialize matrices $\mathbf{S}$, $\mathbf{S}_{N}$,
    $\overline{\mathbf{S}}$, $\overline{\mathbf{S}}_{N}$ to zero
    \FOR{$i\gets 1, \ldots, N_{\psi}$}
    \STATE for $\psi_{i}$ generate integration vectors $\mathbf{q}$,
    $\mathbf{q}_{C}$, $\mathbf{q}_{S}$ using
    Algorithm~\ref{alg:integration}
    \FOR{$j\gets1,\ldots N_{W}$}
    \STATE increment row $j$ of $\mathbf{S}$ by
    $W_{ji}w_{i}^{(\psi)}\mathbf{q}$
    \STATE increment row $j$ of $\mathbf{S}_{N}$ by
    $V_{ji}w_{i}^{(\psi)}\mathbf{q}$
    \STATE increment row $j$ of $\overline{\mathbf{S}}$ by
    $\overline{W}_{ji}w_{i}^{(\psi)}\mathbf{q}+\overline{W}_{ji}^{(C)}\mathbf{q}_{C}
    +\overline{W}_{ji}^{(S)}\mathbf{q}_{S}$
    \STATE increment row $j$ of $\overline{\mathbf{S}}_{N}$ by
    $\overline{V}_{ji}w_{i}^{(\psi)}\mathbf{q}+\overline{V}_{ji}^{(C)}\mathbf{q}_{C}
    +\overline{V}_{ji}^{(S)}\mathbf{q}_{S}$
    \ENDFOR
    \ENDFOR
  \STATE output $N_{W}\times N_{Q}$ matrices $\mathbf{S}$,
  $\overline{\mathbf{S}}$, $\mathbf{S}_{N}$, $\overline{\mathbf{S}}_{N}$
  \end{algorithmic}
\end{algorithm}

\begin{algorithm}
  \caption{Evaluation of radiated field}
  \label{alg:field}
  \begin{algorithmic}
    \STATE input matrices $\mathbf{S}$, $\overline{\mathbf{S}}$,
    $\mathbf{S}_{N}$, $\overline{\mathbf{S}}_{N}$ from
    Algorithm~\ref{alg:shifts}, source terms
    $\boldsymbol{\phi}_{i}=p(t_{i},\theta,\phi)$ and
    $\boldsymbol{\phi}_{i}^{(n)}=\partial \phi(t_{i},\theta,\phi)/\partial n$ at
    $N_{Q}$ nodes on sphere surface
    \STATE initialize acoustic pressure and normal derivative $p$ and
    $p^{(n)}$ to zero
    \FOR{$i=1,\ldots,$}
    \STATE $p_{i,\ldots,i+N_{W}}\gets \mathbf{S}\boldsymbol{\phi}_{i}
    +\mathbf{S}_{N}\boldsymbol{\phi}^{(n)}_{i}$ 
    \STATE $p^{(n)}_{i,\ldots,i+N_{W}}\gets
    \overline{\mathbf{S}}\boldsymbol{\phi}_{i} 
    +\overline{\mathbf{S}}_{N}\boldsymbol{\phi}^{(n)}_{i}$ 
    \ENDFOR
  \end{algorithmic}
\end{algorithm}

\section{Author declarations}
\label{sec:declarations}

The author has no conflicts of interest to declare.

\section{Data availability}
\label{sec:data:availability}

Code implementing the method of the paper and generating the results
presented is available upon request to the author.


\end{document}